\documentclass{article}
\usepackage{amssymb}
\newtheorem{theorem}{Theorem}

\newtheorem{definition}[theorem]{Definition}
\newtheorem{example}[theorem]{Example}
\newtheorem{remark}[theorem]{Remark}

\def\QED{\quad\blackslug\lower 8.5pt\null}

\newcommand{\crazy}[2]{\displaystyle{\mathop{#1}_{#2}}
\vphantom{\displaystyle{#1}}}

\begin{document}

\begin{center}
{\Large\bf MULTIDIMENSIONAL {\mathversion{bold} $(n+1)$}-WEBS }

 \vspace*{3mm}

{\Large\bf WITH REDUCT REDUCIBLE
 SUBWEBS  }

\vspace*{4mm}

{\large VLADISLAV V. GOLDBERG}
\end{center}

{\small{\bf Abstract}.
{\em  We find three characterizations for a
 multidimensional $(n+1)$-web $W$ possessing a reduct
 reducible subweb: its closed form equations, the integrability of
 an invariant distribution associated with $W$,
 and the relations between the components
 of its torsion tensor. In the case of codimension one, the latter
 criterion establishes a relation with solutions of
 a system of nonlinear second-order PDEs.
 Some particular cases of   this system were considered by Goursat in 1899.}

 \vspace*{8mm}

\textbf{1.} Let $W(n+1,n,r)$ be an $(n+1)$-web defined on a
differentiable manifold $X^{nr}$ of dimension $nr$
by $n + 1$ foliations
$\lambda_{\xi},\;\; \xi = 1,...,n+1$ of codimension $r$.
Each foliation $\lambda_{\xi}$
can be defined
by the completely integrable system of Pfaffian
equations
\begin{equation}\label{eq:1}
 \crazy{\omega}{\xi}^{i}=0,\;\; \xi=1,\ldots,n+1;\;\;
 i=1,\ldots,r.
\end{equation}
The 1-forms $\crazy{\omega}{\alpha}^i,\;\;
\alpha = 1, \ldots, n$,
define a co-frame in the tangent bundle $T(X^{nr})$ and satisfy
the following structure equations:
\begin{equation}\label{eq:2}
d\crazy{\omega}{\alpha}^{i} = \crazy{\omega}{\alpha}^{j} \wedge
\omega_{j}^{i} +  \sum_{\beta\neq\alpha}
\crazy{a}{\alpha\beta}^{i}_{jk} \;
\crazy{\omega}{\alpha}^{j} \wedge \crazy{\omega}{\beta}^{k},
\end{equation}
\begin{equation}\label{eq:3}
d\omega_{j}^{i} - \omega_j^k \wedge \omega_k^i =
\sum_{\alpha,\beta = 1}^{n}
\crazy{b}{\alpha\beta}_{jkl}^{i}\;\crazy{\omega}{\alpha}^{k}
\wedge \crazy{\omega}{\beta}^{l} ,
\end{equation}
where the quantities $\crazy{a}{\alpha\beta}^{i}_{jk}$ and
$\crazy{b}{\alpha\beta}_{jkl}^{i}$ are connected by
certain relations (see [G 73] or [G 74] or [G 88], Section \textbf{1.2}).
We indicate some of these relations:
\begin{equation}\label{eq:4}
\nabla\crazy{a}{\alpha\beta}^{i}_{jk} = \sum_{\gamma = 1}^{n}
(\crazy{a}{\alpha\beta\gamma}^{i}_{jkl} +
\crazy{a}{\alpha\beta}^{i}_{mk}
\crazy{a}{\gamma\alpha}^{m}_{lj} +
\crazy{a}{\alpha\beta}^{i}_{jm}
\crazy{a}{\beta\gamma}^{m}_{kl})\,\crazy{\omega}{\gamma}^{l} ,
\;\; \alpha \neq \beta,
\end{equation}
\begin{equation}\label{eq:5}
\crazy{a}{\alpha\beta}^{i}_{jk} =
\crazy{a}{\beta\alpha}^{i}_{kj},
\end{equation}
\begin{equation}\label{eq:6}
\sum_{\alpha ,\beta=1}^{n} \crazy{a}{\alpha\beta}_{jk}^{i} = 0,
\end{equation}
\begin{equation}\label{eq:7}
\crazy{b}{\alpha\beta}^{i}_{jkl} =
\frac{1}{2}(\crazy{a}{\gamma\alpha\beta}_{jkl}^{i} -
\crazy{a}{\beta\gamma\alpha}_{ljk}^{i}),\;\;
\gamma \neq \alpha,\beta,
\end{equation}
where $\nabla\crazy{a}{\alpha\beta}^{i}_{jk} =
d\crazy{a}{\alpha\beta}^{i}_{jk}
- \crazy{a}{\alpha\beta}^{i}_{lk} \,\omega_j^l
- \crazy{a}{\alpha\beta}^{i}_{jl} \,\omega_k^l
+ \crazy{a}{\alpha\beta}^{l}_{jk} \,\omega_l^i.$
The quantities $\crazy{a}{\alpha\beta}^{i}_{jk}$ and
$\crazy{b}{\alpha\beta}_{jkl}^{i}$ form tensor fields in
the tangent bundle
 $T(X^{nr})$ which are called, respectively, the
{\it torsion} and {\it curvature tensors} of the
web $W(n+1,n,r)$.

Equations (7) and (4)
show that if $n > 2$, the curvature tensor
of  $W(n+1,n,r)$ is expressed in terms of the Pfaffian
derivatives of its torsion tensor and the torsion tensor
itself. This is not the case for $n = 2$, because in formulas
(7), $\gamma \neq \alpha,\beta$.

\textbf{2.} In [G 73] (see also [G 74] or [G 88], Section
\textbf{1.5} or [G 93], Section \textbf{2.8})
a {\em reduct  $(k+1)$-subweb} $[n+1, 1,\ldots, k]$ $(k < n)$,
of the web $W$ was defined as a $(k+1)$-web cut on  the intersection
of $n-k$ leaves $F_{\sigma},\; \sigma =
k+1, \ldots, n; \,  1 < k < n$, of the foliations $\lambda_\sigma$
by the leaves $F_{t},\;\; t = 1, \ldots, k, n+1$,
of the other $k+1$ foliations $\lambda_t$.

There are  ${n+1 \choose k+1}$ reduct
$(k+1)$-subwebs of the web $W$.
For simplicity of notation, in this paper
we  consider the reduct $(k+1)$-subweb
indicated above. It is  defined by the system
 of equations
\begin{equation}\label{eq:8}
 \crazy{\omega}{\sigma}^{i}=0, \;\; \sigma = k + 1,  \ldots, n.
\end{equation}

We can see from the structure equations
 (2) that the torsion tensor  of the $(k+1)$-subweb
$[n+1, k+1,\ldots, n]$
is a subtensor $\crazy{a}{uv}^{i}_{jk}, \; u, v = 1, \ldots , k,$
of the torsion tensor $\crazy{a}{\alpha\beta}^{i}_{jk},
\alpha, \beta = 1, \ldots , n,$ of the web $W$.

\textbf{3.} Suppose that in some domain $D \subset  X^{nr}$
a web $W (n + 1, n, r)$ is defined by
the closed form  equations
\begin{equation}\label{eq:9}
x_{n+1}^i = f^i (x_1^{j_1},\ldots,x_n^{j_n}),\;\;i,j_1,\ldots,j_n=1,\ldots,r ,
\end{equation}
where
\begin{equation}\label{eq:10}
\mbox{det}\left(
\frac{\partial{f^i}}{\partial{x_{\alpha}^{j_{\alpha}}}}
\right) \neq 0 ,\;\;\;\;\; \alpha=1,2,\ldots,n.
\end{equation}
In this case the subweb
$[n+1, 1,\ldots, k]$  is defined by the equations
$$
x_\sigma^i = c_\sigma^i, \;\;\sigma = k + 1, \ldots , n,
$$
where $ c_\sigma^i$ are constants.

For simplicity of notation, we  write equations (10) in the
form
\begin{equation}\label{eq:11}
x_{n+1} = f (x_1,\ldots, x_n), \;\; \mbox{det}\left(
\frac{\partial f}{\partial x_{\alpha}},
\right) \neq 0,
\end{equation}
assuming that $x_\xi, \, \xi = 1, \ldots , n,$ and $f$ are vector
functions: $$ x_{\xi} = (x_{\xi}^i), \;\; f = (f^i), \;\;\;\; i =
1, \ldots , r. $$ Then  the subweb
$[n+1, 1,\ldots, k]$  is
defined by the equations
$$
x_\sigma = c_\sigma, \;\;\sigma = k + 1, \ldots , n,
$$
where $ c_\sigma = (c_\sigma^i)$ are constants.

Without loss of generality, in this paper
we  assume the $(k+1)$-subweb $[n+1, 1,\ldots, k]$
is {\em reducible} of type
\begin{equation}\label{eq:12}
x_{n+1} = f (x_1, \ldots, x_k)
= F (x_1, \ldots, x_l, g (x_{l+1}, \ldots , x_k))
\end{equation}
or $(l+1, k)$-{\em reducible} (for the definition of general reducibility see
 [G 76] or [G 88], Section \textbf{4.1}).
 This definition is valid for any $(k+1)$-web.
 However, for the web $[n+1, 1,\ldots, k]$,
 considered as a reduct  $(k+1)$-subweb
 of the web $W (n+1, n, r)$, this definition must be modified as
 follows:
\begin{equation}\label{eq:13}
\begin{array}{ll}
x_{n+1}\!\!\!\! &= f (x_1, \ldots, x_k, c_{k+1}, \ldots ,
c_{n})\\
\!\!\!\!&= F (x_1, \ldots, x_l, g (x_{l+1}, \ldots , x_k, c_{k+1}, \ldots ,  c_{n}),
c_{k+1}, \ldots ,  c_{n}).
\end{array}
\end{equation}

\textbf{4.} Now we will prove the following theorem.
\begin{theorem} Let $W (n+1, n, r)$ be an $(n+1)$-web
given on $nr$-dimensional manifold $X^{nr}$. Then the following
four statements are equivalent:
\begin{description}
\item[\textbf{(i)}] $W (n+1, n, r)$ possesses  the reduct
 $(l+1, k)$-reducible $(k+1)$-subweb $[n+1, 1,\ldots, k]$
 of type $(13)$.
\item[\textbf{(ii)}] The components
 $\crazy{a}{uv}^{i}_{jk}, \; u, v = 1, \ldots , k,$
 of the  torsion tensor of  $W (n+1, n, r)$
 satisfy the conditions
\begin{equation}\label{eq:14}
\crazy{a}{pa}^{i}_{jk} = \crazy{a}{pb}^{i}_{jk}, \;\; p = 1,
\ldots , l; \, a, b = l + 1, \ldots , k.
\end{equation}
\item[\textbf{(iii)}] The distribution defined by the equations
\begin{equation}\label{eq:15}
\crazy{\omega}{\sigma}^{i} = 0, \;\; \crazy{\omega}{l+1}^{i}
+  \crazy{\omega}{k}^{i} = 0, \;\;
\sigma = k +  1,\ldots , n,
\end{equation}
is integrable.
\item[\textbf{(iv)}] The closed form equation of the web  $W (n+1, n, r)$
has the form
\begin{equation}\label{eq:16}
\begin{array}{ll}
x_{n+1}\!\!\!\! &= F (x_1, \ldots, x_k, x_{k+1}, \ldots ,
x_{n})\\
\!\!\!\!&= f (x_1, \ldots, x_l, g (x_{l+1}, \ldots , x_k, x_{k+1}, \ldots ,  x_{n}),
x_{k+1}, \ldots ,  x_{n}),
\end{array}
\end{equation}
where $f$ and $g$ are arbitrary functions of $n-k+l$ and $n-l$
variables, respectively.
\end{description}
\end{theorem}

{\sf Proof.} In fact, the  implications
\textbf{(i)} $\Longleftrightarrow$ \textbf{(iii)} and
 \textbf{(ii)} $\Longleftrightarrow$
\textbf{(iii)} are immediate consequences of  \S 3 of [G 76] (or
Theorems \textbf{4.1.4}, p. 141, and \textbf{1.8.6}, p. 45, of [G 88]),
and the implication  \textbf{(i)} $\Longleftrightarrow$
\textbf{(iv)} is obvious, because the subweb  $[n+1, 1,\ldots, k]$ )
  is $(l+1, k)$-reducible and because it is a subweb of $W (n + 1, n, r)$
  defined by the equations $x_\sigma = c_\sigma$.
\rule{3mm}{3mm}

\textbf{Remark} It is worth noting that the condition
\begin{equation}\label{eq:17}
\begin{array}{ll}
x_{n+1}\!\!\!\! &= f (x_1, \ldots, x_k, x_{k+1}, \ldots ,
x_{n})\\
\!\!\!\!&= f (x_1, \ldots, x_l, g (x_{l+1}, \ldots , x_k),
x_{k+1}, \ldots ,  x_{n})
\end{array}
\end{equation}
is stronger than condition (16):
\begin{itemize}
\item It implies that the entire web
$W (n+1, n, r)$ (not only its subweb $[n+1, 1,\ldots, k]$)
  is $(l+1, k)$-reducible. For this kind of reducibility the
  torsion tensor of $W$ satisfies not  only conditions (14)  but also
the conditions
\begin{equation}\label{eq:18}
\crazy{a}{\sigma a}^{i}_{jk} = \crazy{a}{\sigma b}^{i}_{jk}, \;\;
\sigma = k + 1, \ldots , n; \, a, b = l + 1, \ldots , k.
\end{equation}

\item It also implies that not only the distribution defined by
 equations (15) is integrable but also the larger distribution
defined by the equations
\begin{equation}\label{eq:19}
 \crazy{\omega}{l+1}^i +  \crazy{\omega}{k}^i = 0
\end{equation}
is integrable.
\end{itemize}

This is an immediate consequence of the results of [G 76], \S 3
(or [G 88], Section \textbf{4.1}).

\textbf{5.} If $r = 1$, i.e.,  if a web $W$ is of codimension one,
then $x_\xi, F, f, g$ are scalar functions. In this case
it follows from  [G 76], \S 3 (or [G 88], Section \textbf{4.1}) that conditions
(14) are equivalent to the following system of second-order nonlinear PDEs:
\begin{equation}\label{eq:20}
\frac{F_{pa}}{F_{pb}} = \frac{F_{a}}{F_{b}}, \;\; p = 1,
\ldots , l; \, a, b = l + 1, \ldots , k.
\end{equation}
It follows from our results that the general solution of this system has the form
(16). In one direction this can be verified directly: (16) implies
(20). However, the converse statement (the general solution of
system (20) is given by the functions of form (16)) is not so
obvious. So, {\em closed form equations
of  webs $W$ with the reduct
$(l+1, k)$-reducible $(k+1)$-subweb $[n+1, 1,\ldots, k]$
of type $(13)$ provide this general solution.}

\textbf{6.} Goursat [Go 99] considered the single PDE
which, in our notation, has the form
$$
\frac{F_{12}}{F_{13}} = \frac{F_{2}}{F_{3}}
$$
for  general $n$ and for the particular cases $n = 4, 5$.
This corresponds to the case $l = 1, k = 3$ in equation
(20). The general solution of this equation is
$$
z = f (x_1, g (x_2, x_3, x_4, \ldots, x_n), x_4, \ldots , x_n),
$$
where $f$ and $g$ are arbitrary functions of  $n - 1$
variables each.
For $n = 4, 5$, this solution has the form
$z = f (x_1, g (x_2, x_3, x_4), x_4)$ (where $f$ and $g$
are arbitrary functions of 3 variables each), and
$z = f (x_1, g (x_2, x_3, x_4, x_5), x_4, x_5)$
(where $f$ and $g$
are arbitrary functions of 4 variables each).

Our system (20) and its general solution (16) are more
general. Moreover, these can be generalized for the case
when $r > 1$, i.e., when $x_\xi, F, f, g$ are vector functions.

\noindent {\em Author's address}:

\noindent
Vladislav V. Goldberg\\
Department of Mathematics\\
New Jersey Institute of Technology\\
Newark, N.J. 07102, U.S.A.

\vspace*{2mm}
\noindent
 E-mail address: vlgold@m.njit.edu

\end{document}